\documentclass[english,12pt]{smfart}

\usepackage{smfthm}
\usepackage{smfenum}
\usepackage{amssymb}
\usepackage{amscd}
\usepackage{amsthm}

\newcommand{\on}{\operatorname}

\newcommand{\Qp}{\mathbf{Q}_p}
\newcommand{\Zp}{\mathbf{Z}_p}

\newcommand{\Fp}{\mathbf{F}_p}
\newcommand{\Fpp}{\mathbf{F}_{p^2}}

\newcommand{\OO}{\mathcal{O}}
\newcommand{\MM}{\mathfrak{m}}
\newcommand{\Qbar}{\overline{\mathbf{Q}}}

\newcommand{\Qpbar}{\overline{\mathbf{Q}}_p}

\newcommand{\aplus}{\mathbf{A}^+}
\newcommand{\bcris}{\mathbf{B}_{\mathrm{cris}}}
\newcommand{\bplus}{\mathbf{B}^+}
\newcommand{\dcris}{\mathbf{D}_{\mathrm{cris}}}
\newcommand{\dfont}{\mathbf{D}}
\newcommand{\nwach}{\mathbf{N}}

\newcommand{\ra}{\rightarrow}
\renewcommand{\phi}{\varphi}

\newcommand{\Fil}{\mathrm{Fil}}
\renewcommand{\geq}{\geqslant}
\renewcommand{\leq}{\leqslant} 
\newcommand{\aQp}{\mathbf{A}_{\Qp}}
\newcommand{\bQp}{\mathbf{B}_{\Qp}}

\begin{document}

\title[Families of crystalline representations]
{Construction of some families of $2$-dimensional crystalline representations}

\author{Laurent Berger}
\address{Harvard Dept of Maths  \\
      One Oxford Street \\
      Cambridge, MA 02138 \\ USA}
\email{laurent@math.harvard.edu}
\urladdr{www.math.harvard.edu/\~{}laurent}

\author{Hanfeng Li}
\address{Department of Mathematics\\
University of Toronto\\
Toronto, ON M5S 3G3 \\
Canada}
\email{hli@fields.toronto.edu}

\author{Hui June Zhu}
\address{Department of Mathematics\\
McMaster University\\
Hamilton, ON L8S 4K1 \\ 
Canada}
\email{zhu@cal.berkeley.edu}


\date{October 2003}

\begin{abstract}
We construct explicitly some analytic families of \'etale
$(\varphi,\Gamma)$-modules, which give 
rise to analytic families of $2$-dimensional
crystalline representations. As an 
application of our constructions, we
verify some conjectures of Breuil 
on the reduction modulo $p$ of those
representations, and extend some results (of 
Deligne, Edixhoven, Fontaine and Serre) on the
representations arising from modular forms.
\end{abstract}

\begin{altabstract}
Nous construisons explicitement des familles analytiques de 
$(\varphi,\Gamma)$-modules \'etales, qui donnent lieu \`a des familles
analytiques de repr\'esentations cristallines de dimension $2$. 
Comme application de nos
constructions, nous v\'erifions des conjectures de Breuil
quant \`a la r\'eduction modulo $p$ de ces repr\'esentations, et nous
\'etendons des r\'esultats (de Deligne, Edixhoven, Fontaine et Serre) 
sur les repr\'esentations associ\'ees aux formes modulaires.
\end{altabstract}

\maketitle

\setcounter{tocdepth}{2}
\tableofcontents

\setlength{\baselineskip}{18pt}

\section{Introduction}
Throughout this article $p$ is a prime number and $\Qpbar$ is an
algebraic closure of $\Qp$. We are interested in $p$-adic representations of
$\on{Gal}(\Qpbar/\Qp)$ and we  
use the language of crystalline
representations (see for example \cite{Bu88sst}).

\subsection{Main results}
The purpose of this article is to construct some explicit 
analytic families of $2$-dimensional crystalline representations of
$\on{Gal}(\Qpbar/\Qp)$. More
precisely, let $k \geq 2$ and for $a_p \in \MM_E$ (where $\MM_E$ is the
maximal ideal of the ring of integers of a finite extension $E \subset
\Qpbar$ of $\Qp$)
let $D_{k,a_p}$ be the filtered
$\phi$-module given by $D_{k, a_p} = E e_1 \oplus E e_2$ where:
\[ \begin{cases} \phi(e_1) = p^{k-1} e_2 \\
\phi(e_2) = -e_1 + a_p e_2 
\end{cases}
\qquad\text{and}\qquad
 \Fil^i D_{k, a_p} = \begin{cases}
D_{k, a_p} & \text{if $i \leq 0$,} \\
E e_1 & \text{if $1 \leq i \leq k-1$,} \\
0 & \text{if $i \geq k$.}
\end{cases} \]

By the theorem of Colmez and Fontaine, there exists a crystalline
$E$-linear
representation $V_{k,a_p}$ such that $\dcris(V^*_{k,a_p}) = D_{k,
a_p}$ where $V^*_{k,a_p}=\on{Hom}(V_{k,a_p},E)$ 
is the dual of $V_{k,a_p}$.
The representations $V_{k,a_p}$
are the representations studied by Breuil in \cite{CB02II}. They
are crystalline, irreducible, and their Hodge-Tate weights
are $0$ and $k-1$. One can show 
that if $V$ is any irreducible 
$2$-dimensional crystalline $\Qpbar$-linear
representation,
then there exist $k \geq 2$, $a_p \in \MM_{\Qpbar}$
and $\eta$ a crystalline character such that $V \simeq V_{k,a_p} \otimes
\eta$ (see \cite[prop 3.1.1]{CB02II} for a proof). 

We will give a direct construction of the $V_{k,a_p}$ 
(for those $a_p$ in a small $p$-adic ball around $0$) which
provides a new proof of the theorem of Colmez-Fontaine for these
representations,
along the lines suggested by Fontaine in \cite[B.2.3]{FO91}. We will actually
construct a family of $p$-adic representations, that is a free
$\Qp \otimes_{\Zp} \Zp[[X]]$-module 
\footnote{this definition of a family differs slightly from the one
given in \cite{BC03}}
$\mathbf{V}_k$ of rank $2$ with a continuous
linear action of $G_{\Qp} = \on{Gal}(\Qpbar/\Qp)$ such that if one
sets $m = \lfloor (k-2)/(p-1) \rfloor$, then for any $\alpha \in \MM_{\Qpbar}$ 
we have $\mathbf{V}_k(\alpha) =  V_{k,p^m \alpha}^*$.

The family $\mathbf{V}_k$ gives an explicit $1$-parameter 
deformation of $V_{k,0}^*$ (parameterized by the trace 
of Frobenius). In particular, 
all the representations $V_{k,a_p}$ with
$v_p(a_p) > \lfloor (k-2)/(p-1) \rfloor$ have the same 
semi-simplified reduction modulo $p$ which
partially answers some questions of Breuil 
\cite[conjecture 6.1.1]{CB02II}.

Let $\overline{V}_{k,a_p}$ be the semi-simplification of the reduction
modulo $p$ of a lattice of
$V_{k,a_p}$. Our main result is then (see theorem \ref{reduc}):

\begin{theo}
If $v_p(a_p) > \lfloor (k-2)/(p-1) \rfloor$, 
then $\overline{V}_{k,a_p} \simeq \overline{V}_{k,0}$.
\end{theo}

As a corollary of this result, we can describe the reduction modulo
$p$ of the representations arising from certain modular forms, thus
generalizing some results of Deligne, Edixhoven, Fontaine and Serre.

\subsection{Method of proof}
A $(\phi,\Gamma)$-module is a finite dimensional vector space over a
$2$-dimensional complete local field along with some extra structure: a
Frobenius operator $\phi$ and an 
action of the procyclic group $\Gamma_{\Qp} = 
\on{Gal}(\Qp(\mu_{p^\infty})/\Qp)$, so that a 
$(\phi,\Gamma)$-module is determined by two matrices 
(those of $\phi$ and of a generator of $\Gamma_{\Qp}$)
satisfying some
conditions (the definition of these objects is recalled below in 
\ref{pgmod}). Fontaine has
constructed in \cite[A.3.4]{FO91} an equivalence of categories from the
category of $p$-adic representations to the category of \'etale 
$(\phi,\Gamma)$-modules.

What we will actually construct are explicit families of 
\'etale $(\phi,\Gamma)$-modules. In order to do that, we
only need to construct two matrices satisfying certain properties.
By a theorem of Dee, a family of \'etale
$(\phi,\Gamma)$-modules
gives rise to a family of $p$-adic
representations. 

The question is then: given an \'etale
$(\phi,\Gamma)$-module that corresponds to a $p$-adic
representation $V$, how do we know whether that representation is
crystalline and if it is, how do we compute $\dcris(V)$? The answer to
those questions is given by the theory of Wach modules 
(developed in \cite{NW96, NW97, PC99, LB04} and recalled below in
\ref{pgmod}). 

\subsection{Plan of the article}
In the second chapter, we recall the definitions of
$(\phi,\Gamma)$-modules and of Wach modules. In the third chapter, we
give the construction of the families of Wach modules, and then describe
the families of crystalline representations obtained in this way.
In the fourth chapter we apply those constructions and verify some
conjectures of Breuil. After that, we give an application
to modular forms. In the final chapter, we ask a couple
of open questions related to our constructions.

\vspace{\baselineskip}\noindent\textbf{Acknowledgments}: We thank
C. Breuil, K. Buzzard, G. Chenevier, R. Coleman and P. Colmez for their comments
and encouragements. The research of Zhu is partially supported by an
NSERC discovery grant.

\section{Representations of $G_{\Qp}$}
Let $E$ be a finite extension of $\Qp$, with ring of integers $\OO_E$,
maximal ideal $\MM_E$
and residue field $k_E$. An $E$-linear representation
of $G_{\Qp}$ is a finite dimensional $E$-vector space $V$ with a 
continuous $E$-linear action of $G_{\Qp} = \on{Gal}(\Qpbar/\Qp)$. 
The underlying $\Qp$-vector space of $V$ is then a $\Qp$-linear
representation of $G_{\Qp}$, along with a
$\Qp$-algebra embedding $E \ra \on{End}_{\Qp}(V)$ which commutes with
the action of $G_{\Qp}$. This observation allows us to extend most
constructions on $\Qp$-linear representations to constructions on
$E$-linear representations, and below we recall those 
constructions which we will use.

\subsection{Crystalline representations}
Let $\bcris$ be the ring constructed by Fontaine in
\cite[2.3]{Bu88per}. Recall that if $V$ is a $\Qp$-linear representation,
then $\dcris(V) = (\bcris \otimes_{\Qp} V)^{G_{\Qp}}$ is a filtered
$\phi$-module of dimension $\leq \dim_{\Qp}(V)$. We say that $V$ is
crystalline if equality holds.  Colmez and Fontaine proved in
\cite[th\'eor\`eme A]{CF00} that the functor $V \mapsto \dcris(V)$ is an equivalence
of categories from the category of crystalline representations to the
category of admissible filtered $\phi$-modules.
 
We say that an $E$-linear representation
$V$ is crystalline if and only if the underlying
$\Qp$-linear representation is crystalline. In this case, $\dcris(V)$ is
an $E$-vector space with an $E$-linear Frobenius and a filtration by
$E$-vector spaces: it is an admissible
$E$-linear filtered $\phi$-module and the
theorem of Colmez-Fontaine then implies
that the functor $V \mapsto \dcris(V)$ is an equivalence of categories
from the category of crystalline $E$-linear representations to the
category of admissible $E$-linear filtered $\phi$-modules.

\subsection{Wach modules and $(\phi,\Gamma)$-modules}\label{pgmod}
Let $\Gamma_{\Qp} = \on{Gal}(\Qp(\mu_{p^\infty})/\Qp)$ and let $\chi$
be the cyclotomic character, so that $\chi: \Gamma_{\Qp} \ra \Zp^*$ is
an isomorphism.
Let $\pi$ be a variable and let
$\aQp$ be the ring of power series $\sum_{i=-\infty}^{+\infty} a_i
\pi^i$ such that $a_i \in \Zp$ and $a_i \ra 0$ as $i \ra
- \infty$. This is a local ring with maximal ideal $(p)$ and its
field of fractions is $\bQp=\aQp[1/p]$. Those two rings are endowed
with a Frobenius $\phi$ defined by $\phi(\pi)=(1+\pi)^p-1$ and an
action of $\Gamma_{\Qp}$ defined by
$\gamma(\pi)=(1+\pi)^{\chi(\gamma)}-1$ for $\gamma \in \Gamma_{\Qp}$.

A $(\phi,\Gamma)$-module is an $\aQp$-module of finite rank, with
semi-linear $\phi$ and continuous
action of $\Gamma_{\Qp}$ which commute with each other.
We say that a $(\phi,\Gamma)$-module $D$ is \'etale if $\phi(D)$
generates $D$ over $\aQp$. Recall that Fontaine
has constructed in \cite[A.3.4]{FO91} a functor $T \mapsto \dfont(T)$ which
associates to any $\Zp$-representation $T$ of $G_{\Qp}$ an \'etale 
$(\phi,\Gamma)$-module, and that the functor $T \mapsto \dfont(T)$ is
then an equivalence of categories.
By inverting $p$, one also gets an equivalence of categories between
the category of ($\Qp$-linear)
$p$-adic representations and 
the category of \'etale $(\phi,\Gamma)$-modules over
$\bQp$. 

If $E$ is a finite extension of $\Qp$, 
we extend the Frobenius and the action of $\Gamma_{\Qp}$ to $E
\otimes_{\Qp} \bQp$ by $E$-linearity.
By restricting our attention to $\OO_E$-modules (or $E$-linear
representations), we then get an equivalence of categories from the
category of $\OO_E$-modules (or $E$-linear
representations) to the category of
$(\phi,\Gamma)$-modules over $\OO_E \otimes_{\Zp} \aQp$ 
(or $(\phi,\Gamma)$-modules over $E \otimes_{\Qp} \bQp$), given by $T
\mapsto \dfont(T)$. The inverse functor will be denoted by $\dfont \mapsto
T(\dfont)$ (or $\dfont \mapsto V(\dfont)$).

If $V$ is a crystalline representation, we can be pretty specific 
about what its $(\phi,\Gamma)$-module looks like. Let $\aplus_{\Qp} =
\Zp[[\pi]]$ and $\bplus_{\Qp} = \Qp 
\otimes_{\Zp} \aplus_{\Qp} = \aplus_{\Qp}[1/p]$. 
The following is proved in \cite[\S II.1, \S III.4]{LB04} (for $\Qp$-linear
representations but the $E$-linear case follows at once):

If $V$ is an $E$-linear representation, then $V$ is crystalline 
with Hodge-Tate weights in $[a,b]$ if and
only if there exists an $E \otimes_{\Qp} \bplus_{\Qp}$-module
$\nwach(V)$ contained in $\dfont(V)$ such that:
\begin{enumerate}
\item $\nwach(V)$ is free of rank $d=\dim_E(V)$ over $E \otimes_{\Qp} \bplus_{\Qp}$;
\item The action of $\Gamma_{\Qp}$ preserves $\nwach(V)$ and is
  trivial on $\nwach(V)/\pi \nwach(V)$;
\item $\phi(\pi^b \nwach(V)) \subset \pi^b \nwach(V)$ and $\phi(\pi^b
  \nwach(V)) / \pi^b \nwach(V)$ is killed by $q^{b-a}$ where
  $q=\phi(\pi)/\pi$. 
\end{enumerate}

We say that a crystalline representation is positive if $\on{Fil}^0
\dcris(V) = \dcris(V)$, or equivalently if $b \leq 0$.
If $V$ is a positive crystalline representation 
then we can take $b=0$ above and if
we endow $\nwach(V)$ with the filtration
$\Fil^i \nwach(V) = \{ x \in \nwach(V),\
\phi(x) \in q^i \nwach(V) \}$, then $\nwach(V)/\pi \nwach(V)$ is an
$E$-linear filtered $\phi$-module and by \cite[\S III.4]{LB04} 
we then have an isomorphism: 
$\nwach(V)/\pi \nwach(V) \simeq \dcris(V)$.

If $T$ is a $G_{\Qp}$-stable lattice in $V$, then $\nwach(T) = \dfont(T) \cap
\nwach(V)$  is an $\OO_E \otimes_{\Zp} \aplus_{\Qp}$-lattice in
$\nwach(V)$ 
(note that $\OO_E \otimes_{\Zp} \aplus_{\Qp} = \OO_E[[\pi]]$)
and by \cite[\S III.4]{LB04}
the functor $T \mapsto \nwach(T)$ gives a bijection
between the $G_{\Qp}$-stable lattices $T$ in $V$ and the 
$\OO_E \otimes_{\Zp} \aplus_{\Qp}$-lattices $\nwach(T)$ in $\nwach(V)$
satisfying:
\begin{enumerate}
\item $\nwach(T)$ is free of rank $d=\dim_E(V)$ over $\OO_E \otimes_{\Zp} \aplus_{\Qp}$;
\item The action of $\Gamma_{\Qp}$ preserves $\nwach(T)$;
\item $\phi(\pi^b
  \nwach(T)) \subset \pi^b \nwach(T)$ and $\phi(\pi^b
  \nwach(T)) / \pi^b \nwach(T)$ is killed by $q^{b-a}$. 
\end{enumerate}
Such an object is called a Wach module. The aim of the following
chapter is to construct families of Wach modules.

\section{Families of Wach modules}
In this chapter, we construct $1$-parameter families of Wach
modules. More precisely, fix $k \geq 2$ and define 
$m = \lfloor (k-2)/(p-1) \rfloor$;
we construct a matrix $P(X) \in
\on{M}(2,\Zp[[\pi,X]])$ and for every $\gamma \in
\Gamma_{\Qp}$ a matrix $G_\gamma(X) \in \on{Id} + \pi
\on{M}(2,\Zp[[\pi,X]])$ such that:
\begin{enumerate}
\item For any $\alpha \in \MM_E$, the matrices $P(\alpha)$ and
  $G_{\gamma}(\alpha)$ can be used to define a Wach module
  $\nwach_{k,\alpha}$ which corresponds to a 
  crystalline representation. 
\item If one sets $a_p =\alpha p^m$, then the  
  $2$-dimensional $E$-linear filtered $\phi$-module
  associated to that crystalline representation
  is given by $D_{k, a_p} = E e_1 \oplus E e_2$ 
  where:
\[ \begin{cases} \phi(e_1) = p^{k-1} e_2 \\
\phi(e_2) = -e_1 + a_p e_2 
\end{cases}
\qquad\text{and}\qquad
 \Fil^i D_{k, a_p} = \begin{cases}
D_{k, a_p} & \text{if $i \leq 0$,} \\
E e_1 & \text{if $1 \leq i \leq k-1$,} \\
0 & \text{if $i \geq k$.}
\end{cases} \]
\end{enumerate}
These are the duals of the 
representations defined by Breuil in \cite[\S 3.1]{CB02II}. They
are crystalline, irreducible, and their Hodge-Tate weights
are $0$ and $-(k-1)$.

\subsection{Construction of Wach modules}
Let us now construct this family of Wach modules.
Recall that $q=\phi(\pi)/\pi$. For $n \geq 1$, we define $q_n =
\phi^{n-1}(q)$ so that $q_1 = q$.
Let \[ \lambda_+ = \prod_{n \geq 0} \frac{\phi^{2n+1}(q)}{p} = 
\frac{q_2}{p} \times \frac{q_4}{p} \times \frac{q_6}{p} \times \cdots
\quad\text{and}\quad
\lambda_-  = \prod_{n \geq 0} \frac{\phi^{2n}(q)}{p} =
\frac{q_1}{p} \times \frac{q_3}{p} \times \frac{q_5}{p} \times \cdots
 \]

\begin{prop}\label{holo}
The functions $\lambda_+$ and $\lambda_- \in \Qp[[\pi]]$ satisfy the
following properties:
\begin{enumerate}
\item $\lambda_+(0) = \lambda_-(0) = 1$; \label{pt_1}
\item $\lambda_-/\gamma(\lambda_-)$ and $\lambda_+/\gamma(\lambda_+)
  \in 1+\pi \Zp[[\pi]]$;  \label{pt_2}
\item $\phi(\lambda_-)=\lambda_+$ and $\phi(\lambda_+)=\lambda_- /
  (q/p)$; \label{pt_3}
\item Let $m = \lfloor (k-2)/(p-1) \rfloor$; 
if we write $p^m (\lambda_-/\lambda_+)^{k-1} = \sum_{i \geq 0} z_i
  \pi^i$ and we define $z=z_0 +z_1 \pi + \cdots + z_{k-2} \pi^{k-2}$ 
then $z \in \Zp[[\pi]]$. \label{pt_4}
\end{enumerate}
\end{prop}

\begin{proof}
Since $q_n(0)=p$, point (\ref{pt_1}) is obvious. Point
(\ref{pt_2}) follows from the fact that $\gamma(q)/q \in  1+\pi \Zp[[\pi]]$ and that:
\[ \frac{\lambda_+}{\gamma(\lambda_+)} = \prod_{n \geq 0}
\phi^{2n+1}\left(\frac{q}{\gamma(q)}\right) 
\quad\text{and}\quad
 \frac{\lambda_-}{\gamma(\lambda_-)} = \prod_{n \geq 0}
\phi^{2n}\left(\frac{q}{\gamma(q)}\right). \] 
Point (\ref{pt_3}) follows immediately from the definitions.

Let us now prove (\ref{pt_4}). Let $R$ be the set of
power series $\sum_{i\geq 0} a_i \pi^i$ such that $a_i \in \Qp$ and
$v_p(a_i) + i/(p-1) \geq 0$. One checks easily that $R$ is a ring and
that $q_n/p$ and $p/q_n \in R$ for all $n \geq 1$ so that
$(\lambda_-/\lambda_+)^{k-1} 
\in R$. This implies that if we write $p^m(\lambda_-/\lambda_+)^{k-1} 
= \sum_{i \geq 0} z_i \pi^i$ then $v_p(z_i) + i/(p-1) \geq m$. In
particular, we have $v_p(z_i) \geq m - (k-2)/(p-1)$ if
$i=0,\cdots,k-2$. Since $m = \lfloor (k-2)/(p-1) \rfloor$, 
this implies that $v_p(z_i) > -1$ and so
that $v_p(z_i) \geq 0$. Therefore, $z \in \Zp[[\pi]]$.
\end{proof}

We define \[ P(X) = 
\begin{pmatrix}
0 & -1 \\
q^{k-1} & Xz
\end{pmatrix}
\quad\text{and}\quad
G^{(k-1)}_\gamma = 
\begin{pmatrix}
\left(\frac{\lambda_+}{\gamma(\lambda_+)}\right)^{k-1} & 0 \\
0 & \left(\frac{\lambda_-}{\gamma(\lambda_-)}\right)^{k-1}
\end{pmatrix}. \]

\begin{lemm}\label{start}
We have:
\[ P(X) \phi(G^{(k-1)}_\gamma) -
G^{(k-1)}_\gamma \gamma(P(X)) 
= \begin{pmatrix}
0 & 0 \\ 0 & \pi^{k-1} \star
\end{pmatrix} 
\in \pi^{k-1} \on{M}(2,\Zp[[\pi,X]]). \]
\end{lemm}

\begin{proof}
Note that $\phi$ and $\gamma \in \Gamma_{\Qp}$ act trivially on $X$.
A direct computation using the fact that $\phi(\lambda_+)  = \lambda_-
/ (q/p)$ shows that 
\[  P(X) \phi(G^{(k-1)}_\gamma) = 
\begin{pmatrix} 
0 & -  \left(\frac{\lambda_+}{\gamma(\lambda_+)}\right)^{k-1} \\
\gamma(q)^{k-1}
\left(\frac{\lambda_-}{\gamma(\lambda_-)}\right)^{k-1} & 
Xz  \left(\frac{\lambda_+}{\gamma(\lambda_+)}\right)^{k-1}
\end{pmatrix} \]
and that
\[ G^{(k-1)}_\gamma \gamma(P(X)) = 
\begin{pmatrix} 
0 & - \left(\frac{\lambda_+}{\gamma(\lambda_+)}\right)^{k-1} \\
\gamma(q)^{k-1}
\left(\frac{\lambda_-}{\gamma(\lambda_-)}\right)^{k-1} &
X \gamma(z) \left(\frac{\lambda_-}{\gamma(\lambda_-)}\right)^{k-1}
\end{pmatrix} \]
and to prove the lemma, we need to show that \[ z  
\left(\frac{\lambda_+}{\gamma(\lambda_+)}\right)^{k-1} - \gamma(z) 
 \left(\frac{\lambda_-}{\gamma(\lambda_-)}\right)^{k-1} 
\in \pi^{k-1} \Zp[[\pi]] = \Zp[[\pi]] \cap \pi^{k-1} \Qp[[\pi]]. \]
It is clear that the above series belongs to $\Zp[[\pi]]$. Finally, by
definition of $z$, we have $z - p^m (\lambda_-/\lambda_+)^{k-1} \in \pi^{k-1}
\Qp[[\pi]]$ which proves the second inclusion.
\end{proof}

\begin{prop}\label{lift}
There exists a unique matrix $G_\gamma(X) \in \on{Id} + \pi
\on{M}(2,\Zp[[\pi,X]])$ such that 
$P(X) \phi(G_\gamma(X)) = G_\gamma(X) \gamma(P(X))$.
\end{prop}

\begin{proof}
We will start by proving the uniqueness of a matrix $G_\gamma(X)$
satisfying the above conditions. 
Assume that there exist
$G_\gamma(X)$ and $G'_\gamma(X)$ satisfying the above conditions and set
$H=G'_\gamma(X)G^{-1}_\gamma(X)$. 
A short computation shows that $H \in \on{Id} + \pi
\on{M}(2,\Zp[[\pi,X]])$ and that
$H P(X) = P(X) \phi(H)$. We will show that this implies $H =
\on{Id}$. 
Assume that this is not the case, and
write $H=\on{Id}+H_\ell \pi^\ell + \cdots$ where $H_\ell
\neq 0$ and $P(X) = P_0 + P_1 \pi + \cdots$. 
The facts that $H P(X) = P(X) \phi(H)$ 
and that $\phi(\pi^\ell)/\pi^\ell \equiv p^\ell \mod{\pi}$
imply that
$H_\ell P_0 = p^\ell P_0 H_\ell$ which in turn implies that $P_0$ has
two eigenvalues the quotient of which is $p^\ell$. Since 
\[ P_0 = 
\begin{pmatrix}
0 & -1 \\ p^{k-1} & Xp^m 
\end{pmatrix}, \]
this is impossible.

Let us now show the existence of $G_\gamma(X)$. 
Recall that by lemma \ref{start}, we have:
\[ P(X) \phi(G^{(k-1)}_\gamma) -
G^{(k-1)}_\gamma \gamma(P(X)) 
= \begin{pmatrix}
0 & 0 \\ 0 & \pi^{k-1} \star 
\end{pmatrix} 
\in \pi^{k-1} \on{M}(2,\Zp[[\pi,X]]), \]
and a direct computation shows that
\[ \begin{pmatrix}
0 & 0 \\ 0 & \pi^{k-1} \star 
\end{pmatrix} \gamma(P(X)^{-1}) = 
\frac{1}{q^{k-1}} \begin{pmatrix}
0 & 0 \\ 0 & \pi^{k-1} \star 
\end{pmatrix} 
\begin{pmatrix}
Xz & 1 \\ -q^{k-1} & 0 
\end{pmatrix}
= \begin{pmatrix}
0 & 0 \\ - \pi^{k-1} \star & 0 
\end{pmatrix} \] 
so that:
\[ G^{(k-1)}_\gamma - P(X) \phi(G^{(k-1)}_\gamma)
\gamma(P(X)^{-1}) = \pi^{k-1} R^{(k-1)} \in 
\pi^{k-1} \on{M}(2,\Zp[[\pi,X]]). \]
We will prove by
recurrence on $\ell \geq k-1$ that there exists $G^{(\ell)}_\gamma  = G^{(\ell)}_\gamma(X)$ such
that $G^{(\ell)}_\gamma \equiv G^{(\ell-1)}_\gamma \mod{\pi^{\ell-1}}$ and
\[ G^{(\ell)}_\gamma - P(X) \phi(G^{(\ell)}_\gamma)
\gamma(P(X)^{-1}) = \pi^{\ell} R^{(\ell)} \in \pi^{\ell}
\on{M}(2,\Zp[[\pi,X]]). \]
We already know this for $\ell=k-1$ and for $\ell \geq k$ we need to
find $H^{(\ell)} \in  \on{M}(2,\Zp[[X]])$ such that if we set
$G_\gamma^{(\ell)} = G_\gamma^{(\ell-1)} + \pi^{\ell-1} H^{(\ell)}$ then
\[ G^{(\ell)}_\gamma - P(X) \phi(G^{(\ell)}_\gamma)
\gamma(P(X)^{-1}) \in \pi^{\ell}
\on{M}(2,\Zp[[\pi,X]]). \]

Bearing in mind that $\phi(\pi)=q\pi$,
this gives the following equations for $H^{(\ell)}$:
\[ \pi^{\ell-1} R^{(\ell-1)} + \pi^{\ell-1} H^{(\ell)} - P(X)
\pi^{\ell-1} \phi(H^{(\ell)}) q^{\ell-1} \gamma(P(X)^{-1}) 
\in  \pi^{\ell} \on{M}(2,\Zp[[\pi,X]]). \]
Notice that since $\ell \geq k$, we have $q^{\ell-1} \gamma(P(X)^{-1})
\in  \on{M}(2,\Zp[[\pi,X]])$.
Therefore, we can divide the above
equation by $\pi^{\ell-1}$ and we see that we only need to find
$H^{(\ell)} \in  \on{M}(2,\Zp[[X]])$ such that 
\[ H^{(\ell)} - P_0 H^{(\ell)} (p^{\ell-1} P_0^{-1}) \equiv - R^{(\ell-1)} \mod{\pi}. \]
Notice that since $\ell \geq k$, we have $p^{\ell-1} P_0^{-1} \in
p^{\ell-k} \on{M}(2,\Zp[[X]])$.
This implies that the operator $H \mapsto H - P_0 H (p^{\ell-1} P_0^{-1})$
is a bijection on $\on{M}(2,\Zp[[X]])$ if $\ell \geq k+1$ and to
finish the proof, we only need to show that this also holds if $\ell=k$.
However, a direct computation shows that modulo $(p,X)$, this operator
becomes:
\[\begin{pmatrix}
h_{11} & h_{12} \\ h_{21} & h_{22} 
\end{pmatrix} 
\mapsto \begin{pmatrix}
h_{11} & h_{12} + h_{21} \\ h_{21} & h_{22} 
\end{pmatrix}, \] which is obviously invertible.
This shows that at each step, we can ``lift'' $G_\gamma^{(\ell-1)}$ to
one and only one $G_\gamma^{(\ell)}$. To finish the proof, we take
$G_\gamma(X) = \lim_{\ell \ra + \infty} G^{(\ell)}_\gamma(X)$.
\end{proof}

\subsection{Families of crystalline representations}
First, we check that the above construction does give rise to a Wach
module, and therefore to a crystalline representation. If $\alpha
\in \MM_E$, then the matrices $P(\alpha)$ and $G_\gamma(\alpha)$ belong to 
$\OO_E[[\pi]]$.

\begin{prop}\label{eval}
If $\alpha \in \MM_E$ and
$\gamma$, $\eta \in \Gamma_K$ then $G_{\gamma \eta}(\alpha) = G_\gamma(\alpha)
\gamma(G_\eta(\alpha))$ and 
\[ P(\alpha) \phi(G_\gamma(\alpha)) = G_\gamma(\alpha)
\gamma(P(\alpha)) \] so that one can use the matrices $P(\alpha)$ and
$G_\gamma(\alpha)$ to define
a Wach module $\nwach_{k,\alpha}$ over $\OO_E[[\pi]]$.
\end{prop} 

\begin{rema}\label{same}
When $\alpha=0$, the Wach module $\nwach_{k,0}$ coincides with the
last example of \cite[Appendice A]{LB04}.
\end{rema}

\begin{proof}
We already know that $P(X) \phi(G_\gamma(X)) = G_\gamma(X)
\gamma(P(X))$ and if $\gamma$, $\eta \in \Gamma_K$ 
then $G_{\gamma \eta}(X)$ and $G'_{\gamma \eta}(X) = G_\gamma(X)
\gamma(G_\eta(X))$ both satisfy the conditions of proposition
\ref{lift} so that they are equal. We then define $\nwach_{k,\alpha}$
as the free $\OO_E[[\pi]]$-module of rank $2$ with basis $n_1$, $n_2$
as follows:
$\nwach_{k,\alpha} = \OO_E[[\pi]] n_1 \oplus \OO_E[[\pi]] n_2$.
We then endow it with a Frobenius $\phi$ and an action of $\Gamma_{\Qp}$ by
deciding that the matrix of $\phi$ with respect to the basis
$(n_1,n_2)$ is
$P(\alpha)$ and that
the matrix of $\gamma \in \Gamma_{\Qp}$ is
$G_\gamma(\alpha)$. 
\end{proof}

\begin{defi}\label{tandv}
If $\alpha \in \MM_E$, define $a_p = p^m \alpha$.
Let $V_{k,a_p}$ be the crystalline $E$-linear representation
such that $E \otimes_{\OO_E} \nwach_{k,\alpha} = \nwach(V_{k,a_p}^*)$ and let 
$T_{k,a_p}$ be the $\OO_E$-lattice in $V_{k,a_p}$ such that 
$\nwach_{k,\alpha} = \nwach(T_{k,a_p}^*)$.
\end{defi}

\begin{prop}\label{visdkap}
The filtered $\phi$-module $E \otimes_{\OO_E}
(\nwach_{k,\alpha} / \pi \nwach_{k,\alpha})$ is isomorphic to
the filtered $\phi$-module $D_{k,a_p}$ described at the beginning 
of the chapter, so that we have $\dcris(V_{k,a_p}^*) = D_{k,a_p}$.
\end{prop}

\begin{proof}
Recall that $D_{k, a_p} = E e_1 \oplus E e_2$ 
  where:
\[ \begin{cases} \phi(e_1) = p^{k-1} e_2 \\
\phi(e_2) = -e_1 + a_p e_2 
\end{cases}
\qquad\text{and}\qquad
\Fil^i D_{k,a_p} = \begin{cases}
D_{k,a_p} & \text{if $i \leq 0$,} \\
E e_1 & \text{if $1 \leq i \leq k-1$,} \\
0 & \text{if $i \geq k$.}
\end{cases} \] and that \[ 
P(\alpha) = \begin{pmatrix} 0 & -1 \\ q^{k-1} & \alpha z \end{pmatrix}
\quad\text{so that}\quad
(P(\alpha) \mod{\pi}) = \begin{pmatrix} 0 & -1 \\ p^{k-1} & a_p 
\end{pmatrix}. \]
We will show that the map $n_i \mapsto e_i$ gives the required
isomorphism from $E \otimes_{\OO_E}
(\nwach_{k,\alpha} / \pi \nwach_{k,\alpha})$ to $D_{k,a_p}$. It is
obvious from the above that this is an isomorphism of $\phi$-modules,
and to finish the proof we need to compute the filtration on
$\nwach_{k,\alpha}$. The proposition will then follow from:
\[ \Fil^i \nwach_{k,\alpha} = \begin{cases}
\nwach_{k,\alpha} & \text{if $i \leq 0$,} \\
n_1 \OO_E[[\pi]] \oplus \pi^i n_2  \OO_E[[\pi]] & \text{if $1 \leq i \leq k-1$,} \\
\pi^{i-(k-1)}  n_1 \OO_E[[\pi]]  \oplus \pi^i n_2  \OO_E[[\pi]] & \text{if $i \geq k$.}
\end{cases} \]
Recall that $\Fil^i \nwach = \{ x \in \nwach,\
\phi(x) \in q^i \nwach \}$. It is obvious that $\Fil^i
\nwach_{k,\alpha} = \nwach_{k,\alpha}$ if $i \leq 0$.

Next, choose $1 \leq i \leq k-1$ and write $x = x_1 n_1 + x_2 n_2$. We
have $\phi(x) = -\phi(x_2) n_1 + (q^{k-1} \phi(x_1) + \alpha z \phi(x_2))
n_2$ and therefore, $x \in \Fil^i \nwach_{k,\alpha}$ if and only if
$q^i | \phi(x_2)$ which is equivalent to $\pi^i | x_2$. 
If $i \geq k$, the proof is completely similar.
\end{proof}

\begin{rema}
By Dee's theorem \cite[theorem 2.1.27]{JD01} 
(see also \cite{BC03}), an \'etale family of $(\phi,\Gamma)$-modules (in our
case, $\nwach_{k,X}$ over $\Zp[[X]]$) gives rise to an analytic family of $p$-adic
representations (in our case, a free $\Zp[[X]]$-module $\mathbf{T}_k$ of rank $2$
with a continuous linear action of $G_{\Qp}$).

The family $\mathbf{V}_k = \Qp \otimes_{\Zp} \mathbf{T}_k$
is a $1$-dimensional subspace (parameterized by $a_p$, the trace of
Frobenius) of the space of all
representations deforming $V_{k,0}^*$.
\end{rema}

\section{Applications}
\subsection{Reduction modulo $p$}
We will now use the constructions of the previous chapter to verify
some conjectures of Breuil
about the reduction of the representations
$V_{k,a_p}$ defined above. Recall that $a_p = p^m \alpha$ where $m =
\lfloor (k-2)/(p-1) \rfloor$ and $\alpha \in \MM_{\Qpbar}$.

\begin{theo}\label{reduc}
Given $\alpha \in \MM_E$, the two $k_E$-representations 
$k_E \otimes_{\OO_E} T_{k,a_p}$ and $k_E \otimes_{\OO_E}
T_{k,0}$ are isomorphic.
\end{theo}

\begin{proof}
Given $\alpha \in \MM_E$, we have $G_\gamma(\alpha) \equiv
G_\gamma(0) \mod{\MM_E}$ and $P(\alpha) \equiv P(0) \mod{\MM_E}$
because $P$, $G_\gamma \in \on{M}(2,\aplus_{\Qp}[[X]])$. 
If $\dfont_{k,\alpha} = \aQp \otimes_{\aplus_{\Qp}} \nwach_{k,\alpha}$ 
is the $(\phi,\Gamma)$-module associated to $\nwach_{k,\alpha}$ then we
have \[ \on{Hom}_{\phi,\Gamma}\left(\dfont_{k,0},k_E((\pi))\right) =
\on{Hom}_{\phi,\Gamma}\left(\dfont_{k,\alpha},k_E((\pi))\right). \]
Since $k_E \otimes_{\OO_E} T_{k,a_p} =
T(\on{Hom}_{\phi,\Gamma}(\dfont_{k,\alpha},k_E((\pi))))$ 
where $T$ is
Fontaine's functor (since $T(\cdot)$ is an exact functor), 
we see that  $k_E \otimes_{\OO_E} 
T_{k,a_p} \simeq k_E \otimes_{\OO_E}
T_{k,0}$.
\end{proof}

\begin{rema}\label{remreduc}
\begin{enumerate}
\item The constant $m=\lfloor (k-2)/(p-1) \rfloor$ is not necessarily
  the ``best'' possible one for every $k$. 
  Indeed, in proposition \ref{holo}-\ref{pt_4} one
  only needs to take $m$ such that $(p^m (\lambda_-/\lambda_+)^{k-1}
  \mod{\pi^{k-1}})$ has coefficients in $\Zp$. For example, if $k=p+1$,
  then one may take $m=0$; \label{hl}
\item The same argument will show that if $i \geq 1$ and 
$\alpha_1$, $\alpha_2 \in \MM_E$
satisfy $v_E(\alpha_1-\alpha_2) \geq i$ 
and $a_p^{(j)}=p^m \alpha_j$ then $T_{k,a_p^{(1)}} \equiv
T_{k,a_p^{(2)}} \mod{\MM_E^i}$;
\item Instead of using the $(\phi,\Gamma)$-modules, one could
  also use the fact that by Dee's theorem, the $T_{k,a_p}^*$'s come from
  evaluating $\mathbf{T}_k$ at $\alpha$ where $a_p = p^m \alpha$ and the
  above congruences follow at once.
\end{enumerate}
\end{rema}

We can now prove some congruences which were conjectured
(with slightly sharper bounds on $v_p(a_p)$) by Breuil (See Breuil's
\cite[conj 6.1.1]{CB02II}). 

If $a_p \in \MM_E$ and $v_p(a_p) > m = \lfloor 
(k-2)/(p-1) \rfloor$, then there exists $\alpha \in \MM_E$ such that
$a_p = \alpha p^m$ and the following is then 
a direct consequence of theorem \ref{reduc} (and remark
\ref{remreduc}-\ref{hl}
for $k=p+1$)
applied to $m=0,1,2$
(for $k\leq p$, the first congruence is also a direct consequence of
``Fontaine-Laffaille''). 
\begin{coro}
The following congruences hold:
\begin{enumerate}
\item If $k \leq p+1$, then
$\overline{V}_{k,a_p} \simeq \overline{V}_{k,0}$ if $v_p(a_p) > 0$;
\item If $k \leq 2p-1$, then
$\overline{V}_{k,a_p} \simeq \overline{V}_{k,0}$ if $v_p(a_p) > 1$;
\item If $k \leq 3p-2$, then 
$\overline{V}_{k,a_p} \simeq \overline{V}_{k,0}$ if $v_p(a_p) > 2$.
\end{enumerate}
\end{coro}

To keep this article reasonably self-contained, 
let us recall the explicit description of
$\overline{V}_{k,0}$ given by Breuil 
in \cite[prop 6.1.2]{CB02II}. Recall that
$\chi$ is the cyclotomic character; $\omega_2$ stands for the
fundamental character of level $2$, $\mu_{\pm\sqrt{-1}}$ is the
unramified character of $G_{\Qp}$ which sends $\on{Frob}_p$ to
$\pm\sqrt{-1}$ and $\on{ind}$ is the induction from $\Fpp$
to $\Fp$.

\begin{prop}\label{ap0}
We have:
\begin{enumerate}
\item If $(p+1) \nmid (k-1)$, then $\overline{V}_{k,0}=\on{ind}(\omega_2^{k-1})$;
\item If $(p+1) | (k-1)$, then \[ \overline{V}_{k,0} = \begin{pmatrix}
\mu_{\sqrt{-1}} & 0 \\ 0 & \mu_{-\sqrt{-1}}  \end{pmatrix} \otimes
  \chi^{(k-1)/(p+1)}. \] 
\end{enumerate}
\end{prop}

This proposition is itself a consequence of \cite[prop 3.1.2]{CB02II}
which describes $V_{k,0}$. 

\subsection{Application to modular forms}
In this paragraph, we give an application of theorem \ref{reduc} above
to the reduction modulo $p$ of the representations attached to certain
modular forms. This way, we can generalize results of Deligne,
Edixhoven, Fontaine and Serre (see \cite[\S 2]{BE92} and
\cite[th\'eor\`eme 6.2.3]{CB02II}). 
We follow \cite[\S 6.2]{CB02II} very closely.

Fix an embedding $\Qbar
\hookrightarrow \Qpbar$.
Let $f$ be a normalized cuspidal
modular form over $\Gamma_1(N)$, 
of level $N \geq 1$ coprime to $p$, 
of character $\eta$ and
of weight $k \geq 2$.
We also assume that
$f$ is a newform (hence an eigenform) 
for the Hecke operators
$T_\ell$ for all primes $\ell \nmid N$. 
The elements $a_\ell(f) \in \OO_{\Qpbar}$
are defined by $T_\ell f = a_\ell(f) f$ 
and we assume in this
paragraph that $a_p(f) \in \MM_{\Qpbar}$ 
(in other words, that $f$ is not ordinary).

Let $V_f$ be the restriction to $G_{\Qp}$ of the
$p$-adic representation attached to $f$ and let $\overline{V}_f$
denote the semi-simplification of its reduction modulo $p$. 
Recall that we have
\[ V_f \simeq V_{k,a_p(f) \eta^{1/2}(p)} 
\otimes \eta^{1/2}, \]
where $\eta^{1/2}$
is an unramified character of $G_{\Qp}$ 
whose square is $\eta$ (see 
\cite[th\'eor\`eme 6.2.1]{CB02II} for a statement and
\cite[\S 4.3]{BM02} for precise
references).  
We can therefore apply theorem \ref{reduc} and
get:

\begin{theo}\label{modform}
Under the preceding hypotheses, we have $\overline{V}_f \simeq
\overline{V}_{k,0} \otimes 
\eta^{1/2}$ if $v_p(a_p) > \lfloor
(k-2)/(p-1) \rfloor$. 
\end{theo}

This way, we get an explicit description of $\overline{V}_f$
using \cite[prop 6.1.2]{CB02II} recalled above in proposition 
\ref{ap0}.

\section{Open questions}
In this last chapter, we give a few open questions related to the
construction of $(\phi,\Gamma)$-modules. 

\begin{enumerate}
\item It is clear that our methods should extend to more general
  situations: given a filtered $\phi$-module, can one find a Wach
  module whose reduction modulo $\pi$ is that filtered $\phi$-module?
  This article gives some evidence of how one can do that, but in
  general one problem is the following: given a lattice $T$ in a crystalline
  representation $V$, then $\nwach(T) / \pi \nwach(T)$ is a lattice in
  $\dcris(V)$. Which lattices do we get in this way? Do we get all
  strongly divisible lattices (see \cite{FL82} for a definition) 
  in $\dcris(V)$ this way? When the length of the filtration is $\leq
  p-1$, the answer is yes (see \cite[\S 3.2]{NW97} and \cite[\S
  V.2]{LB04}). 
\item In the ramified case or in the semi-stable case, we know even
  less. However, when the length of the filtration is $\leq p-1$ (or
  $\leq p-2$ in the semi-stable case) those representations are
  described by ``Breuil modules'' (see \cite{CB98,CB00}). 
  How does one compute the
  $(\phi,\Gamma)$-module of such a representation directly from its
  Breuil module? In the unramified crystalline case, 
  when the length of the
  filtration is $\leq p-1$, then a Breuil module is equivalent to the
  data of a strongly divisible lattice 
  and the computation 
  of the associated Wach module is done
  in \cite[\S 3.2]{NW97} and in \cite[V.2]{LB04}.
\end{enumerate}

\end{document}